\documentclass[reqno]{amsart}

\usepackage{amssymb, calrsfs, graphics, graphicx, enumerate, enumitem, url, xcolor, hyperref, mathrsfs, listings, comment}
\usepackage{tabularx}
\usepackage[ruled, lined, linesnumbered, longend]{algorithm2e}
\usepackage[justification=centering]{caption}

\newtheorem{theorem}{Theorem}[section]
\newtheorem{lemma}[theorem]{Lemma}

\newtheorem{corollary}[theorem]{Corollary}

\numberwithin{equation}{section}

\theoremstyle{remark}
\newtheorem*{remark}{Remark}

\lstdefinestyle{CStyle}{
    basicstyle=\footnotesize,
    breakatwhitespace=false,         
    breaklines=true,                 
    captionpos=b,                    
    keepspaces=true,                 
    numbers=left,                    
    numbersep=5pt,                  
    showspaces=false,                
    showstringspaces=false,
    showtabs=false,                  
    tabsize=2,
    language=C
}
\title[Burgess' bounds]{Burgess' bounds  with prime moduli}

\author[B. Kerr]{Bryce Kerr}
\address{School of Science\\ The University of New South Wales Canberra, Australia}
\thanks{Supported by Australian Research Council Discovery Project DP230100534.}
\email{bryce.kerr@unsw.edu.au}

\author[J. Punch]{James Punch}
\address{School of Science\\ The University of New South Wales Canberra, Australia}
\email{j.punch@unsw.edu.au}

\author[T.~S. Trudgian]{Timothy S. Trudgian}
\address{School of Science\\ The University of New South Wales Canberra, Australia}
\thanks{Supported by Australian Research Council Discovery Project DP240100186.}
\email{t.trudgian@unsw.edu.au}

\begin{document}

\begin{abstract}
Let $S_{\chi}(N)= \sum_{n\leq N} \chi(n)$, where $\chi$ is a Dirichlet character modulo a prime $p$. For integral $r\geq 2$, Burgess gave bounds of the shape $S_{\chi}(N) \ll N^{1-1/r} p^{(r+1)/(4r^{2})} \log p$. The factor of $\log p$ was reduced by several authors. We show that it can be removed completely.
\end{abstract}
\maketitle
\section{Introduction and statement of result}
Let $\chi$ be a Dirichlet character modulo a prime $p$, whence either $\chi(n) = 0$ or $|\chi(n)| =1$. A central problem in analytic number theory is to give bounds on the sum $S_{\chi}(N) = \sum_{n\leq N} \chi(n)$. Trivially, we have that $|S_{\chi}(N)| \leq N$, which is useful for small $N$. The P\'{o}lya--Vinogradov bound of $|S_{\chi}(N)| \leq \sqrt{p}\log p$ is useful for large $N$. Bounds for `intermediate' $N$ were obtained by Burgess \cite{Burgess62}.

Indeed, Burgess \cite{Burgess63} proved  
\begin{equation}\label{gumnut}
    S_{\chi}(N) \ll N^{1-1/r} p^{(r+1)/(4r^{2})} f(p),
    \end{equation}
for $f(p) = \log p$. While the upper bound on $S_{\chi}(N)$ is expected to be $\ll_{\varepsilon} \sqrt{N}q^{\varepsilon}$ (see \cite{McGownTrevino} for more details), the only progress made since Burgess' work is in reducing the term $f(p)$ in (\ref{gumnut}). To that end, Friedlander and Iwaniec \cite[p.\ 366]{Fried} showed that taking $f(p)=(\log p)^{1/2r}$ is possible. Explicit bounds appear in Iwaniec and Kowalski \cite{IK} (with $f(p) = (\log p)^{1/r}$), and in Booker \cite{Booker} (only for quadratic characters $\chi$, with $f(p) = (\log p)^{1/2r}$). These have since been improved by McGown \cite{Macca}, Trevi\~{n}o \cite{Trev}, and Francis \cite{Forro}. Further improvements came with $f(p) = (\log p)^{1/4r}$ in Kerr et al. \cite{Kerr}. While these (and subsequent improvements on $f(p)$) have not been made explicit, we note that Burgess' bounds where $\chi$ has composite modulus have recently been made explicit by Bordignon \cite{Bordo}, Jain-Sharma et al.\ \cite{Sharmo}, Hasanalizade et al.\ \cite{Marto}, and the second author \cite{Puncho}.

The bound on $f(p)$ was further reduced to $f(p) = (\log p)^{(\theta+\varepsilon)/r}$ for $\theta =0.0832\ldots$ and any $\varepsilon>0$ by de la Bret\`eche et al. \cite{Regis}. The main focal point of \cite{Regis} is an optimal result for G\'{a}l sums. We pursue a different, and more elementary, line of argument in this article.  Since we expect the bound in (\ref{gumnut}) to be far from the truth, one may certainly hope to eliminate\footnote{ 
While this has not been done before our work, we note that via an alternative approach, in some applications (see, e.g., \cite{KevTim}) logarithmic powers have been removed completely.} the logarithms in $f(p)$ completely. This is what we do in our main result below.

\begin{theorem}\label{th1}
Let $\chi$ be a non-principal Dirichlet character modulo a prime $p$, let $r\geq 2$ be an integer, and let $\alpha>1$ be some real number. Suppose that $p>r^{\frac{2r}{2r-1}}$ and that $N\leq p^{\frac{1}{2}+\frac{1}{4r}}$. Let $$C(r)=\frac{\alpha\cdot (4r-2)^\frac{r}{r-1}\cdot \frac{4r^{1/2r}}{\left(1-2^{\frac{2r}{2r-1}-2}\right)^{1-1/2r}}}{(4r-2)^\frac{2r-1}{2r-2}-2(4r-2)^\frac{1}{2r-2}},$$ and suppose that \begin{align}\label{condition}\frac{C(r)^r}{(4r-2)^\frac{r}{r-1}}\cdot p^\frac{r-3}{4r}\cdot \left(\frac{(\alpha -1)p^{1/2r}}{\alpha}-\frac{1}{r}\right)\geq r.\end{align} Then for each integer $r\geq 2$, we have
$$\sum_{n\leq N} \chi(n) \leq C(r) N^{1-1/r} p^{(r+1)/(4r^{2})}.$$
\end{theorem}
\begin{remark}
    It is evident that if \eqref{condition} holds for $p=p_0$, then it must hold for $p\geq p_0$. For example, when $r=2$ and $\alpha = 1.002$, then \eqref{condition} holds for $p\geq 3\times 10^{12}$ and we have $C(2)<37$.
\end{remark}
\begin{corollary} \label{cor}
    Let $\chi$ be a non-principal Dirichlet character modulo a prime $p$ and suppose that $N\leq p^{\frac{5}{8}}$. Then $$\sum_{n\leq N}\chi(n)\leq 37 N^{1/2}p^{3/16}.$$
\end{corollary}
\begin{proof}
    For $p\geq 3\times 10^{12}$, this follows from Theorem \ref{th1}. Also, if $37N^{1/2}p^{3/16}\geq N$, this just follows from the trivial bound. Since $N\leq p^{5/8}$, this happens if $$N^{1/2}\leq 37 p^{3/16}\impliedby p^{5/16}\leq 37p^{3/16} \impliedby p\leq 37^{8}=(3.5\ldots)\times 10^{12}.$$
\end{proof}
We note briefly that the point of our article is to show that the $\log p$ factors can be removed in the Burgess bound via Theorem \ref{th1}. That result and Corollary \ref{cor} are indeed explicit, but are unlikely to be useful in applications in their current form. For example, Corollary \ref{cor} only improves on the bounds in Table 3 in \cite{Forro}, when $p>10^{8046}$.  

The outline of this article is as follows.
We introduce some preliminary results in \S \ref{jam}. These allow us to prove Theorem \ref{th1} and Corollary \ref{cor} in \S \ref{macgyver}. Finally, in \S \ref{jambon} we outline some potential improvements to our result.

\section{Preliminary lemmas}\label{jam}
Let $A,N$ be chosen such that $A<N$ and $AN<p$. For $x \in \mathbb{F}_{p}^{*}$, consider
$$\nu(x) = \# \{(n,a):1\leq n \leq N, 1\leq a \leq A, \textrm{ and } ax \equiv n \pmod{p}\}.$$
Since each value of $a$ gives rise to at most one value of $n$, we must have $\nu(x)\leq A$. We show that much more is true in our Lemma \ref{lem} below. First, though, we need to record a preliminary result, which is Lemma 6.2 from \cite{Bagger}.


\begin{lemma}[Bagger et al.]\label{mole}
Let $p$ be prime and $U$ and $V$ be positive integers satisfying $UV < p$. Let $\lambda$ be such that $p\nmid \lambda$, and let $\Delta\geq 1$. Suppose that the congruence
$$v \equiv \lambda u \pmod{p}, \quad 1\leq u \leq U, \quad 1 \leq v \leq V,$$
has at least $\Delta$ solutions. Then there exist positive integers $u_{0}, v_{0}$ such that
$$v_{0} \equiv \lambda u_{0} \pmod{p}, \quad 1 \leq u_{0} \leq \frac{U}{\Delta},  \quad 1 \leq v_{0} \leq \frac{V}{\Delta}, \quad (u_{0}, v_{0}) = 1.$$
\end{lemma}
Using Lemma \ref{mole}, we  improve on the bound $\nu(x) \leq A$ in the following result.

\begin{lemma}\label{lem}
Suppose that $A<N$ and $AN<p$. If $1<\alpha <2$, then $$\sum_{x} \nu(x)^{\alpha} \leq AN\cdot \frac{2^\alpha}{1-2^{\alpha-2}}.$$
\end{lemma}
\begin{proof}
We partition the range of $a$ dyadically. If we write $A\sim N$ to mean $A\in [N, 2N)$, we get
$$\sum_{x} \nu(x)^{\alpha} \leq \sum_{2^{j} \leq A} 2^{\alpha (j+1)} \cdot\#\{x\in \mathbb{F}_{p}^{*} \textrm{ such that }  \nu(x) \sim 2^{j}\}.$$
We now apply Lemma \ref{mole}. 
If $x$ is such that $\nu(x)\sim 2^{j}$, then there exist $a_{0}, n_{0}$ for which
$$ 1 \leq a_{0} \leq A/2^{j}, \quad 1 \leq n_{0} \leq N/2^{j}, \quad n_{0} \equiv x a_{0} \pmod{p}.$$ The pair $(a_0,n_0)$ determines $x$, so there cannot be more than $\frac{AN}{2^{2j}}$ options for $x$ that have $\nu(x)\sim 2^j$. Hence
$$\sum_{x} \nu(x)^{\alpha} \leq  2^\alpha \sum_{2^{j} \leq A} 2^{\alpha j}\frac{AN}{2^{2j}} \leq 2^\alpha AN \sum_{j\geq 0} 2^{(\alpha-2)j}=AN\cdot \frac{2^{\alpha}}{1-2^{\alpha-2}}$$
as required.
\end{proof}
\section{Log-free Burgess}\label{macgyver}

The proof of Theorem \ref{th1}, which we now give, closely follows that of Theorem 12.6 in \cite{IK}.
We want to prove that $\lvert S_{\chi}(N) \rvert \leq C(r)N^{1-1/r}p^{\frac{r+1}{4r^2}}$. This follows from the trivial bound if $N\geq C(r)^r p^{\frac{r+1}{4r}}$. We thus assume that $$C(r)^r p^{\frac{r+1}{4r}}\leq N \leq p^{\frac{1}{2}+\frac{1}{4r}},$$ and induct on $N$. Shifting by positive integers $h\in [1,H]$ yields $$S_{\chi}(N)=\sum_{n\leq N}\chi(n+h)+2\theta E(H),$$ with $E(H)=C(r)H^{1-1/r}p^{\frac{r+1}{4r^2}}$ and $\lvert \theta \rvert\leq 1$. Use shifts $h=ab$ with $1\leq a \leq A$ and $1 \leq b \leq B$; averaging gives $\lvert S_{\chi}(N)\rvert \leq H^{-1}V+2E(H),$ where $H=AB$ and $$V=\sum_{x \pmod{p}}\nu(x)\left \lvert \sum_{1\leq b \leq B}\chi(x+b)\right\rvert.$$ By H\"{o}lder's inequality, we have 

\begin{equation}\label{tbone}
V\leq \left(\sum_{x} v(x)^{2r/(2r-1)}\right)^{1- 1/2r} \left( \sum_{x} \left| \sum_{1\leq b \leq B} \chi (x + b)\right|^{2r} \right)^{1/2r}.
\end{equation}
We treat the two sums over $x$ using Lemma \ref{lem} and \cite[Lemma 12.8]{IK}. To do so, we must assume that $A<N$, $AN<p$ and $B<p$; we validate these assumptions momentarily. We arrive at $$V\leq \left(AN \cdot \frac{2^{\frac{2r}{2r-1}}}{1-2^{\frac{2r}{2r-1}-2}} \right)^\frac{2r-1}{2r}\cdot \left((2rB)^r p+2rB^{2r} \sqrt{p}\right)^\frac{1}{2r}.$$ Select $B=\left\lfloor rp^{1/2r}\right\rfloor$ and $A=\left\lfloor\frac{N}{Krp^{1/2r}}\right\rfloor$ for some real $K>1$. Clearly, $A<N$; further, since $p>r^{\frac{2r}{2r-1}}$, we get $B<p$. Since $N\leq p^{\frac{1}{2}+\frac{1}{4r}}$, we also have $AN<p$.

Then \begin{align*}
    V&\leq N^{1-1/2r} \left(\frac{N}{Krp^{1/2r}}\right)^{1-1/2r}\frac{2}{\left(1-2^{\frac{2r}{2r-1}-2}\right)^{1-1/2r}}\left((2r)^{2r}p^{3/2}\right)^{1/2r} \\
    &=N^{2-1/r}p^{\frac{r+1}{4r^2}}K^{-1+1/2r}\frac{4r^{1/2r}}{\left(1-2^{\frac{2r}{2r-1}-2}\right)^{1-1/2r}}=N^{2-1/r}p^{\frac{r+1}{4r^2}}K^{-1+1/2r}d(r),
\end{align*} say. Clearly, $AB\leq N/K$. Due to the assumed lower bound on $N$, for any $L>K$ we will have $AB\geq N/L$ for sufficiently large $p$. Assume we have such an $L$; we will specify it shortly. It follows that
\begin{align*}
    \lvert S_{\chi}(N)\rvert &\leq \frac{L}{N}N^{2-1/r}p^{\frac{r+1}{4r^2}}K^{-1+1/2r}d(r)+2C(r)\left(\frac{N}{K}\right)^{1-1/r}p^\frac{r+1}{4r^2}\\
    &=N^{1-1/r}p^{\frac{r+1}{4r^2}}\left[\frac{Ld(r)}{K^{1-1/2r}}+\frac{2C(r)}{K^{1-1/r}}\right].
\end{align*}

We will be done if the term in brackets equals $C(r)$; this happens if we choose $$C(r)=\frac{L d(r)}{K^{1-1/2r}-2K^{1/2r}}=\frac{L}{K}\cdot d(r)\cdot \frac{K}{K^{1-1/2r}-2K^{1/2r}}.$$ The final term is minimised when $K=(4r-2)^\frac{r}{r-1}$; we use this $K$ and select $L=\alpha K$. This gives $$C(r)=\frac{\alpha\cdot (4r-2)^\frac{r}{r-1}\cdot d(r)}{(4r-2)^{\frac{2r-1}{2r-2}}-2(4r-2)^\frac{1}{2r-2}}.$$ It remains to establish that $AB\geq N/L$.

This occurs if \begin{align*}
    \left(\frac{N}{Krp^{1/2r}}-1\right)\left(rp^{1/2r}-1\right) &\geq \frac{N}{\alpha K}\\
    \impliedby \frac{N}{K}\left(\frac{\alpha -1}{\alpha}-\frac{1}{rp^{1/2r}}\right)&\geq rp^{1/2r} \\
    \impliedby \frac{C(r)^r}{K} \left(\frac{\alpha -1}{\alpha}\cdot p^{\frac{r-1}{4r}}-\frac{p^{\frac{r-3}{4r}}}{r}\right) &\geq r,
\end{align*} which is implied by \eqref{condition}.

\section{Potential improvements}\label{jambon}

There are several places where our argument could be improved, or optimised further. Since there are some obstacles to these improvements, we present them here as future research directions.

First, in our Lemma \ref{lem} we split the range of $a$ up dyadically. One instead could consider values of $\nu(x)\in[\lambda^{j}, \lambda^{j+1})$, where $\lambda>1$ is some parameter, and then optimise over $\lambda$. The improvement here would be small; for instance, when $r=2$ the constant $37$ would improve to about $36.7$.

Second, one could consider a further refinement of Lemma \ref{lem}. Since $AN < p$, the lattice points under consideration lie on a common line through the origin. If there are at least  $m$ points then there is a primitive point in the rectangle with sides $A/m$ and $N/m$. One could use partial summation to count these primitive points. 

Third, we can, when $r=2$, estimate the fourth moment in (\ref{tbone}) directly. This follows a usual Weil-type argument where one separately estimates exceptional tuples.

Fourth, instead of averaging over positive values of $b$, average over negative values as well, and attenuate by a diminishing weight $\omega_{b}$. We expect the ensuing optimisation problem (of choosing the best functions $\omega_{b}$) to give a substantial saving.

We end with a remark concerning current results in number theory. One expects that the Burgess bounds, with or without the logarithms, are far from the truth. Indeed one expects that $S_{\chi}(N)\ll N^{1/2} q^{\epsilon}$ --- see, e.g., \cite[(12.54)]{IK}. A recent OpenAI result \cite{OAI} shows a sort of quasi-Riemann hypothesis that has bearing on this bound. Consider Dirichlet $L$-functions, $L(s, \chi)$, where $\chi$ is a character mod $q$. The result in \cite{OAI} is that there are no zeroes of $L(s, \chi)$ when $\Re(s) > 7/8$. This would show that $S_{\chi}(N)\ll N^{7/8} q^{\epsilon}$. It remains to be seen what explicit constants could be obtained with this approach. Until then, though, we expect that our Theorem \ref{th1} and Corollary \ref{cor} will remain useful in applications.



\section*{Acknowledgements}
We are grateful to Igor Shparlinski for discussions on this topic over the years. The first and third authors are also grateful to Banff International Research Station, and the organisers of the conference ``Arithmetic, $L$-functions, and Pseudorandomness'' (26w5546), at which we discussed this with several participants.


\begin{thebibliography}{10}

\bibitem{Booker}
A.~R. Booker.
\newblock Quadratic class numbers and character sums.
\newblock {\em Math. Comp.}, 75(255):1481--1492, 2006.

\bibitem{Bordo}
M. Bordignon.
\newblock Partial {G}aussian sums and the {P}\'olya--{V}inogradov inequality for
  primitive characters.
\newblock {\em Rev. Mat. Iberoam.}, 38(4):1101--1127, 2022.

\bibitem{Burgess62}
D.~A. Burgess.
\newblock On character sums and primitive roots.
\newblock {\em Proc. London Math. Soc. (3)}, 12:179--192, 1962.

\bibitem{Burgess63}
D.~A. Burgess.
\newblock On character sums and {$L$}-series. {II}.
\newblock {\em Proc. London Math. Soc. (3)}, 13:524--536, 1963.

\bibitem{Regis}
R. de~la Bret\`eche, M. Munsch, and G.Tenenbaum.
\newblock Small {G}\'al sums and applications.
\newblock {\em J. Lond. Math. Soc. (2)}, 103(1):336--352, 2021.

\bibitem{Forro}
F.~J. Francis.
\newblock An investigation into explicit versions of {B}urgess' bound.
\newblock {\em J. Number Theory}, 228:87--107, 2021.

\bibitem{Fried}
J.~Friedlander and H.~Iwaniec.
\newblock Estimates for character sums.
\newblock {\em Proc. Amer. Math. Soc.}, 119(2):365--372, 1993.

\bibitem{Bagger}
G. K. Bagger, A. R. Booker, B. Kerr, K. McGown, V. Starichkova, and T. Trudgian.
\newblock The determination of norm-Euclidean cyclic cubic fields.
\newblock {\em Preprint available at arXiv:2507.05862}, 2026.

\bibitem{Marto}
E.~Hasanalizade, H.~Lin, G.~Martin, A.~L. Mart\'{\i}nez, and E.~Trevi\~{n}o.
\newblock Explicit Burgess inequalities for cubefree moduli.
\newblock {\em Preprint available at arXiv:2511.17778}, 2025.

\bibitem{IK}
H. Iwaniec and E. Kowalski.
\newblock {\em Analytic Number Theory}, volume~53 of {\em American Mathematical
  Society Colloquium Publications}.
\newblock American Mathematical Society, Providence, RI, 2004.

\bibitem{Sharmo}
N. Jain-Sharma, T. Khale, and M. Liu.
\newblock Explicit {B}urgess bound for composite moduli.
\newblock {\em Int. J. Number Theory}, 17(10):2207--2219, 2021.

\bibitem{Kerr}
B. Kerr, I.~E. Shparlinski, and K.-H. Yau.
\newblock A refinement of the {B}urgess bound for character sums.
\newblock {\em Michigan Math. J.}, 69(2):227--240, 2020.

\bibitem{McGownTrevino}
K. McGown and E. Trevi\~no.
\newblock The least quadratic non-residue.
\newblock In {\em Mexican Mathematicians in the World---Trends and Recent
  Contributions}, volume 775 of {\em Contemp. Math.}, pages 205--231. Amer.
  Math. Soc., Providence, RI, 2021.

\bibitem{Macca}
K.~J. McGown.
\newblock On the constant in {B}urgess' bound for the number of consecutive
  residues or non-residues.
\newblock {\em Funct. Approx. Comment. Math.}, 46:273--284, 2012.

\bibitem{KevTim}
K.~J. McGown and T. Trudgian.
\newblock Explicit upper bounds on the least primitive root.
\newblock {\em Proc. Amer. Math. Soc.}, 148(3):1049--1061, 2020.

\bibitem{OAI}
OpenAI.
\newblock The quasi-Riemann hypothesis: A zero-free half-plane $\Re s > 7/8$.
\newblock {\em Preprint}, 2026.

\bibitem{Puncho}
J. Punch.
\newblock Improved explicit Burgess-strength bounds for composite moduli.
\newblock {\em Preprint available at arXiv:2610.22939}, 2026.

\bibitem{Trev}
E. Trevi\~no.
\newblock The {B}urgess inequality and the least {$k$}th power non-residue.
\newblock {\em Int. J. Number Theory}, 11(5):1653--1678, 2015.

\end{thebibliography}
\end{document}